\documentclass[12pt,reqno]{amsart}

\usepackage[mathlines]{lineno}

\usepackage{times}
\usepackage{amsfonts}
\usepackage{amssymb}
\usepackage{latexsym}
\usepackage{xcolor}

\newtheorem{theorem}{Theorem}

\begin{document}

\author{A. R. Mirotin}
\address{Department of Mathematics and Programming Technologies, Francisk Skorina Gomel State University, Gomel, 246019, Belarus $\&$ Regional Mathematical Center, Southern Federal
University, Rostov-on-Don, 344090, Russia}
\email{amirotin@yandex.ru}

\

\begin{center}
\Large{Some remarks on the solution of  the cell growth equation}
\end{center}

\

\begin{center}
A. R. Mirotin
\end{center}
\

\textsc{Abstract.} The analytical solution to the  initial--boundary value problem for  the cell growth equation was given in the paper  Zaidi A. A., Van Brunt B., Wake G.C., Solutions to an advanced functional partial differential equation of the pantograph type, Proc. R. Soc. A 471: 20140947 ( 2015). In this note, we simplify the arguments given in the paper mentioned above by using the theory of operator semigroups.

\

\section{Introduction}

 The analytical solution to the  initial--boundary value problem for  the cell growth equation was given  at the first time  by  A.A. Zaidi,  B. Van Brunt, and G.C. Wake in \cite{AVBW}. They proved also the uniqueness of the solution and find its asymptotics (for the asymptotic result  see also \cite{PerthameRyzhik}). 
 
A  cell growth model under consideration was developed in \cite{HW}.  See \cite{AVBW} for the history of the issue and a detailed bibliography.

 Let $n(x,t)$  denote
the number density functions of cells of size $x$ at time $t$. Then 

 \begin{equation}\label{CG1}
\frac{\partial n(x,t)}{\partial  t}+g\frac{\partial n(x,t)}{\partial  x}= b\alpha^2n(\alpha x,t)-(b+\mu) n(x,t),
\end{equation}
where $g>0$ is the rate of growth, $\mu>0$ is the rate of death, and $b>0$ is the rate at
which cells divide into $\alpha>1$  equally sized daughter cells.

The above equation is supplemented by a given initial distribution
\begin{equation}\label{IC1}
n(x, 0) = n_0(x), 
\end{equation}
where $n_0$ is a probability distribution function, and the boundary condition,
\begin{equation}\label{BC1}
n(0, t) = 0. 
\end{equation}
The main goal of this note, is to simplify the arguments given in \cite{AVBW} on the existence and uniqueness of the solution of  the  initial--boundary value problem (\ref{CG1}), (\ref{IC1}), (\ref{BC1}) and to discuss its asymptotics  from the point of view  of the theory of   operator  semigroups.

\section{Existence and uniqueness of the solution}
Following  \cite{AVBW} we put
 \begin{equation}\label{CHANGE} 
 n(x, t) = e^{-(b+\mu)t}u(x, t).
\end{equation}

Then,
\begin{equation}\label{CG2}
\frac{\partial u(x,t)}{\partial  t}=-g\frac{\partial u(x,t)}{\partial  x}+ b\alpha^2 u(\alpha x,t)
\end{equation}
 and conditions (\ref{IC1}) and (\ref{BC1}) take the form
 \begin{equation}\label{IC2}
 u(x,0)=n_0(x),
 \end{equation}
 
\begin{equation}\label{BC2}
u(0, t) = 0. 
\end{equation} 
 
 Note that each solution $u$ to  the  initial--boundary value problem (\ref{CG2}), (\ref{IC2}), (\ref{BC2})  for $x\ge 0$  can be extended to a solution to  the initial problem (\ref{CG2}), (\ref{IC2}) for $x\in \mathbb{R}$ if we put $u(x,t)=0$ and $n_0(x)=0$ for all $x\le 0$, $t\ge 0$. Thus, we shall seek
 the solutions $u(x,t)$ to the problem (\ref{CG2}), (\ref{IC2}) for $x\in \mathbb{R}$, $t\ge 0$ putting $n_0(x)=0$ for $x\le 0$.
 
 We rewrite the Cauchy problem (\ref{CG2}) , (\ref{IC2}) with $x\in \mathbb{R}$ in an abstract form in the usual manner. Let  $X=L^p(\mathbb{R})$, $1\le p<\infty$, or  $X=C_{\mathrm{ub}}(\mathbb{R})$ the space of uniformly continuous bounded functions on $\mathbb{R}$ endowed with the $\sup$ norm.
 Consider the vector-valued function $u(t):=u(\cdot,t)$, $u:\mathbb{R}_+\to X$ and the following  operator on $X$:
$$
   \mathcal{H}f(x)= b\alpha^2 f(\alpha x).  
$$
Then the equation (\ref{CG2})  takes the form
$$
\frac{du(t)}{dt}=\left(-g\frac{d}{d  x}+\mathcal{H}\right)u(t).
$$
The operator
$$
 G:=-g\frac{d}{d  x}+\mathcal{H}
$$
 is a generator of a $C_0$-group $S(t)$
   on $X$, since the operator $A=-g d/dx$ with an appropriate domain $D(A)$ \footnote{For instance for $X=L^p(\mathbb{R})$ we  have  $D(A) = \{f \in L^p(\mathbb{R}) : f  \mbox{ absolutely   continuous   and }
   f' \in L^p(\mathbb{R})\}$  (see, e.g., \cite[p. 66]{EN}).}
 is a generator
of a  $C_0$-group of shifts $S_0(t)f(x)=f(x-gt)$ on $X$, $\|S_0(t)\|=1$, and $\mathcal{H}$ is bounded on $X$ \cite[Theorem 13.2.2]{HPh}.
Therefore,
\begin{equation}\label{SOL}
u(t)=S(t)u_0
\end{equation}
is  the unique solution to  the Cauchy problem  
\begin{equation}\label{CPinX}
\frac{du(t)}{dt}=Gu(t),  u(0)=u_0
\end{equation}
for any $u_0\in D(A)$. The formula (\ref{SOL}) gives also a so-called mild solution to  the Cauchy problem   for any $u_0\in X$. 

 We shall assume that $u_0(x)=n_0(x)$ for $ x\in \mathbb{R}_+$ and $u_0(x)=0$ for $x<0$.

The group  $S(t)$
 can be calculated via the Dyson--Phillips series
\begin{equation}\label{DPh}
S(t)=\sum_{n=0}^\infty S_n(t)
\end{equation}
where $\|S_n(t)\|\le \|\mathcal{H}\|^nt^n/n!$ (see 
\cite[(13.2.5)]{HPh}) and 
 \begin{equation}\label{Sn+1}
 S_{n+1}(t)=\int_0^t S_0(t-s)\mathcal{H}S_n(s)ds,\quad n\in \mathbb{Z}_+
 \end{equation}
(see \cite[(13.2.4)]{HPh}, or \cite[Theorem III.1.10]{EN}). 

Thus,
 \begin{equation}\label{SOL2}
 u(x,t)=\sum_{n=0}^\infty S_n(t)u_0(x).
\end{equation}
 Note that $S_0(t)u_0(x)=u_0(x-gt)\ge 0$ for all  $x\in \mathbb{R}$, $t\ge 0$ if $u_0(x)\ge 0$ for all  $x\in \mathbb{R}$, and $S_0(t)u_0(x)=0$  if $x\le 0$, $t\ge 0$. Now it follows from (\ref{Sn+1}) by induction that $u(x,t)$  is non-negative  for all $x\in \mathbb{R}$, $t\ge 0$ and equals to zero for all $x\le 0$, $t\ge 0$.
 
 Moreover, since $\|S_0(t)\|_{X\to X}=1$, it follows (\cite[Corollary of the Theorem 13.2.1]{HPh}, or \cite[Theorem III.1.3]{EN}), that
\begin{equation}\label{EST1}
\|S(t)\|_{X\to X}\le e^{t\|\mathcal{H}\|_{X\to X}}\quad (t\ge 0).
\end{equation}
This yields
\begin{equation}\label{EST}
\|u(\cdot,t)\|_{X}\le e^{t\|\mathcal{H}\|_{X\to X}}\|u_0\|_{X} \quad (t\ge 0).
\end{equation}
In particular, if we assume as in \cite{AVBW} that $u_0\in L^1(\mathbb{R}_+)$ we get for $X=L^1(\mathbb{R})$ that
$$
\|u(\cdot,t)\|_{L^1}\le e^{b\alpha t}\|u_0\|_{L^1}\quad (t\ge 0).
$$
This estimate is consistent with the asymptotics for $u(x,t)$ proven in \cite{AVBW}.

On the other hand, let  $X=C_{\mathrm{ub}}(\mathbb{R})$ and $u_0\in C_{\mathrm{ub}}(\mathbb{R})$. Then we deduce from (\ref{EST}) that 
$$
|u(x,t)|\le e^{b\alpha^2 t}\sup_{\mathbb{R}_+}|u_0| \mbox{ for all } x,t\in \mathbb{R}_+.
$$

In summary, we have the following result.

\begin{theorem}
Let  $X=L^p(\mathbb{R})$, $1\le p<\infty$, or  $X=C_{\mathrm{ub}}(\mathbb{R})$. Let $n_0\in X$ be non-negative. Then formula (\ref{SOL2})
presents a non-negative solution $u$ to  the  initial--boundary value problem (\ref{CG2}), (\ref{IC2}), (\ref{BC2}) such that $u(\cdot,t)\in X$ for $t\ge 0$.
 Moreover, this solution is unique and the estimate (\ref{EST}) holds.
\end{theorem}

{\bf Remark 1.} It follows from (\ref{SOL2}) that
\begin{equation*}
 u(\cdot,t)\approx \sum_{k=0}^n S_k(t)u_0,
\end{equation*}
and by \cite[(13.2.6)]{HPh}
\begin{equation*}
 \|u(\cdot,t)-\sum_{k=0}^n S_k(t)u_0\|_X\le \|\mathcal{H}\|_{X\to X}^{n+1}t^{n+1} \frac{e^{t\|\mathcal{H}\|_{X\to X}}}{(n+1)!}.
\end{equation*}

\section{On the asymptotics of the solution as $t\to\infty$}

It was proven in \cite{AVBW} (cf. \cite{PerthameRyzhik})  for the case $n_0\in L^1(\mathbb{R}_+)$ that 
$u(x,t)\sim e^{b\alpha t}y(x)$  pointwise  as $t\to\infty$.

On can derive several complements to this result from a general theory of operator semigroups, as well. Recall that a function $f$ in $L^1_{loc}(\mathbb{R}_+, X)$
converges to an element $y\in X$ ($X$ is a Banach space) in a sense of Ces\` aro as $t\to\infty$ if
$$
C\mbox{-}\lim_{t\to\infty}f(t):=\lim_{t\to\infty}\frac{1}{t}\int_0^t f(s)ds=y
$$
(convergence in the norm of $X$; see, e.g., \cite{Arendt}). Let  $X=L^1(\mathbb{R})$. Note that  $T(t):=e^{-b\alpha t}S(t)$ is a bounded $C_0$-semigroup in  $L^1(\mathbb{R})$ by (\ref{EST1}) with the generator $B:=G-b\alpha I$ that satisfies the condition (4.5) from \cite[P. 261]{Arendt}.  Let
$u_0\in \mathrm{Ker} B+\overline{\mathrm{Ran} B}$ with $\mathrm{Ker} B$
and $\mathrm{Ran} B$ denoting the kernel and range of $B$. Since $u(\cdot,t)=S(t)u_0$, we have by \cite[Proposition  4.3.1]{Arendt} that $u(\cdot,t)\sim e^{b\alpha t}y$  for some $y\in L^1(\mathbb{R})$  as $t\to\infty$
in a sense that
$$
C\mbox{-}\lim_{t\to\infty} e^{-b\alpha t}u(\cdot,t)=y
$$
(convergence in $L^1(\mathbb{R})$). 
Moreover,  by this Proposition 
$$
y=\lim_{\lambda\downarrow 0}\lambda R(\lambda,G-b\alpha I)u_0=\lim_{\lambda\downarrow 0}\lambda R(\lambda+b\alpha,G)u_0.
$$
Thus, the problem of the asymptotics of the solution of our equation reduces to the  asymptotics of the resolvent of the operator $G$.

The results detailed above  follow mutatis mutandis for several another choices of  $X$ and semigroup generators and
for a more general choice of constants.




\begin{thebibliography}{99}
\bibitem{Arendt}
W. Arendt et.al.,
Vector-valued Laplace Transforms and
Cauchy Problems. Second Edition,  Springer, Basel AG (2011). 

\bibitem{EN}
Engel  K.-J., Nagel R. , One-parameter Semigroups for Linear Evolution Equations, Springer, NY (2000).


\bibitem{HW}
Hall A.J. ,  Wake  G.C., A functional differential equation arising in the
modelling of cell-growth, J. Aust. Math. Soc. Ser. B, 30, 424--435  (1989).



\bibitem{HPh}
E. Hille,  R. S. Phillips, Functional Analysis and Semi-Groups. (revised edition), Amer. Math. Soc. 
Colloq. Publ. Vol. 31, Providence, R. I. (1957).

\bibitem{PerthameRyzhik}
 Perthame B, Ryzhik L. 2005 Exponential decay for the fragmentation or cell-division equation.
J. Diff. Eq. 210, 155--177. (doi:10.1016/j.jde.2004.10.018)


\bibitem{AVBW}
 Zaidi A.A., Van Brunt B., Wake G.C., Solutions to an advanced functional
partial differential equation of the pantograph
type, Proc. R. Soc. A 471: 20140947 ( 2015).

\end{thebibliography}
\end{document}